\documentclass[10pt, a4paper]{article}
\usepackage{amscd}
\usepackage{amsfonts}

\begin{document}

\newtheorem{theorem}{Theorem}
\newtheorem{corollary}[theorem]{Corollary}
\newtheorem{definition}[theorem]{Definition}
\newtheorem{lemma}[theorem]{Lemma}
\newtheorem{proposition}[theorem]{Proposition}
\newtheorem{remark}[theorem]{Remark}
\newtheorem{example}[theorem]{Example}
\newtheorem{notation}[theorem]{Notation}
\def\Qed{\hfill\raisebox{.6ex}{\framebox[2.5mm]{}}\\[.15in]}

\title{A note on Todorov surfaces}


\author{Carlos Rito}

\date{}
\pagestyle{myheadings}
\maketitle

\begin{abstract}

Let $S$ be a {\em Todorov surface}, {\it i.e.}, a minimal smooth surface of general type with
$q=0$ and $p_g=1$ having an involution $i$ such that
$S/i$ is birational to a $K3$ surface and such that the bicanonical map of $S$ is composed with $i.$

The main result of this paper is that, if $P$ is the minimal smooth model of $S/i,$
then $P$ is the minimal desingularization of a double cover of $\mathbb P^2$ ramified over two cubics.
Furthermore it is also shown that, given a Todorov surface $S$, it is possible to construct
Todorov surfaces $S_j$ with $K^2=1,\ldots,K_S^2-1$ and such that $P$ is also the smooth
minimal model of $S_j/i_j,$ where $i_j$ is the involution of $S_j.$
Some examples are also given, namely an example different from the examples presented by Todorov in \cite{To2}.

\noindent 2000 Mathematics Classification: 14J29, 14J28.
\end{abstract}

\section{Introduction}

An {\em involution} of a surface $S$ is an
automorphism of $S$ of order 2. We say that a map is {\em composed with
an involution} $i$ of $S$ if it factors through the double cover $S\rightarrow S/i.$
Involutions appear in many contexts in the study of algebraic surfaces.
For instance in most cases the bicanonical map of a surface of general type
is non-birational only if it is composed with an involution.

Assume that $S$ is a smooth minimal surface of general type with $q=0$ and $p_g\ne 0$
having bicanonical map $\phi_2$ composed with an involution $i$ of $S$ such that
$S/i$ is non-ruled. 
Then, according to \cite[Theorem 3]{Xi}, $p_g(S)=1,$ $K_S^2\leq 8$
and $S/i$ is birational to a $K3$ surface
(Theorem 3 of \cite{Xi} contains the assumption ${\rm \deg}(\phi_2)=2,$ but the result
is still valid assuming only that $\phi_2$ is composed with an involution).

Todorov (\cite{To2}) was the first to give examples of such surfaces.
His construction is as follows.
Consider a Kummer surface $Q$ in $\mathbb P^3,$ {\it i.e.}, a quartic having as only
singularities 16 nodes $a_i.$ The double cover of $Q$ ramified over the intersection of
$Q$ with a general quadric and over the $16$ nodes of $Q$ is a surface of general type
with $q=0,$ $p_g=1$ and $K^2=8.$
Then, choose $a_1,\ldots,a_6$ in general position and let $G$ be the intersection
of $Q$ with a general quadric through $j$ of the nodes $a_1,\ldots,a_6.$
The double cover of $Q$ ramified over $Q\bigcap G$ and over the remaining $16-j$ nodes
of $Q$ is a surface of general type with $q=0,$ $p_g=1$ and $K^2=8-j.$

Imposing the passage of the branch curve by a 7-th node,
one can obtain a surface with $K^2=p_g=1$ and $q=0.$ 
This is the so-called {\em Kunev surface}.
Todorov (\cite{To1}) has shown that the Kunev surface is a bidouble cover
of $\mathbb P^2$ ramified over two cubics and a line.

I refer to \cite{Mo} for an explicit description
of the moduli spaces of Todorov surfaces.

We call {\em Todorov surfaces} smooth surfaces $S$ of general type
with $p_g=1$ and $q=0$ having bicanonical map composed with an involution $i$ of $S$
such that $S/i$ is birational to a $K3$ surface.

In this paper we prove the following:

\begin{theorem}\label{propositionK3}
Let $S$ be a Todorov surface with involution $i$ and $P$ be the smooth minimal model of $S/i.$
Then:
\begin{description}
  \item[a)] there exists a generically finite degree $2$ morphism $P\rightarrow\mathbb P^2$ ramified over two cubics;
  \item[b)] for each $j\in\{1,\ldots,K_S^2-1\},$ there is a Todorov surface $S_j,$
  with involution $i_j,$ such that $K_{S_j}^2=j$ and $P$ is the smooth minimal model of $S_j/{i_j}.$

\end{description}
\end{theorem}

The idea of the proof is the following. First we verify that the evenness of the branch
locus $B'+\sum A_i\subset P$ implies that each nodal curve $A_j$ can only be contained
in a Dynkin graph $\textsf G$ of type $\textsf A_{2n+1}$ or $\textsf D_n.$ Then we use
a Saint-Donat result to show that $A_j$ can be chosen such that the linear system
$|B'-\textsf G|$ is free. This implies b). Finally we conclude that there is a free linear
system $|B_0'|$ with $B_0'^2=2,$ which gives a).

\bigskip
\noindent{\bf Notation and conventions}

We work over the complex numbers;
all varieties are assumed to be projective algebraic.
For a projective smooth surface $S,$ the {\em canonical class} is denoted by
$K,$ the {\em geometric genus} by $p_g:=h^0(S,\mathcal O_S(K)),$ the {\em irregularity}
by $q:=h^1(S,\mathcal O_S(K))$ and the {\em Euler characteristic} by
$\chi=\chi(\mathcal O_S)=1+p_g-q.$

A {\em $(-2)$-curve} or {\em nodal curve} on a surface is a curve
isomorphic to $\mathbb P^1$ such that $C^2=-2$. We say that a curve
singularity is {\em negligible} if it is either a double point or a
triple point which resolves to at most a double point after one
blow-up.

The rest of the notation is standard in algebraic geometry.

\bigskip
\noindent{\bf Acknowledgements}

The author is a collaborator of the Center for Mathematical
Analysis, Geometry and Dynamical Systems of Instituto Superior T\'
ecnico, Universidade T\' ecnica de Lisboa, and
is a member of the Mathematics Department of the
Universidade de Tr\'as-os-Montes e Alto Douro.
This research was partially supported by FCT (Portugal) through
Project POCTI/MAT/44068/2002.

\section{Preliminaries}\label{SaintDonat}
The next result follows from
 \cite[(4.1), Theorem 5.2, Propositions 5.6 and 5.7]{St}.
\begin{theorem}\label{SD}{\rm (\cite{St})}
Let $|D|$ be a complete linear system on a smooth $K3$ surface $F,$ without
fixed components and such that $D^2\geq 4.$ Denote by $\varphi_D$
the map given by $|D|.$ If $\varphi_D$ is non-birational and the surface
$\varphi_D(F)$ is singular then there exists an elliptic pencil $|E|$ such that
$ED=2$ and one of these cases occur:
\begin{description}
    \item [(i)] $D=O_F(4E+2\Gamma)$ where $\Gamma$ is a smooth rational
    irreducible curve such that $\Gamma E=1.$ In this case
    $\varphi_D(F)$ is a cone over a rational normal twisted quartic
    in $\mathbb P^4;$
    \item [(ii)] $D=O_F(3E+2\Gamma_0+\Gamma_1),$ where $\Gamma_0$ and $\Gamma_1$ are smooth rational
    irreducible curves such that $\Gamma_0 E=1,$ $\Gamma_1 E=0$ and
    $\Gamma_0\Gamma_1=1.$ In this case
    $\varphi_D(F)$ is a cone over a rational normal twisted cubic
    in $\mathbb P^3;$
    \item [(iii)]\begin{description}
        \item [a)] $D=O_F(2E+\Gamma_0+\Gamma_1),$ where $\Gamma_0$ and $\Gamma_1$ are smooth rational
    irreducible curves such that $\Gamma_0 E=\Gamma_1 E=1$ and
    $\Gamma_0\Gamma_1=0;$
        \item [b)] $D=O_F(2E+\Delta),$ with $\Delta=2\Gamma_0+\cdots
        +2\Gamma_N+\Gamma_{N+1}+\Gamma_{N+2}$ $(N\geq 0)$, where the
        curves $\Gamma_i$ are irreducible rational curves as in Figure \ref{SDConfig}.

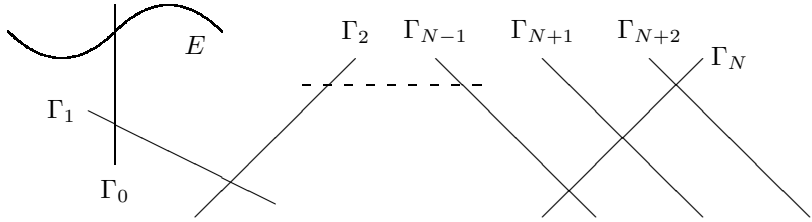
\begin{figure}[htpb]
\begin{center}
\begin{picture}(240,80)(0,0)
\put(-10,30){\line(2,-1){70}}
\put(30,-10){\line(1,1){60}}
\multiput(70,40)(8,0){9}{\line(1,0){3}}
\put(120,50){\line(1,-1){60}}
\put(160,50){\line(1,-1){60}}
\put(200,50){\line(1,-1){60}}
\put(160,-10){\line(1,1){60}}
\put(0,10){\line(0,1){60}}
\put(90,60){\makebox(0,0){$\Gamma_2$}}
\put(120,60){\makebox(0,0){$\Gamma_{N-1}$}}
\put(160,60){\makebox(0,0){$\Gamma_{N+1}$}}
\put(200,60){\makebox(0,0){$\Gamma_{N+2}$}}
\put(230,50){\makebox(0,0){$\Gamma_N$}}
\put(-20,30){\makebox(0,0){$\Gamma_1$}}
\put(0,0){\makebox(0,0){$\Gamma_0$}}
\put(30,55){\makebox(0,0){$E$}}
\qbezier(-40,60)(-20,40)(0,60)
\qbezier(40,60)(20,80)(0,60)
\end{picture}
\end{center}
\caption{Configuration {\bf (iii) b)}}
\label{SDConfig}
\end{figure}
    \end{description}

    In both cases $\varphi_D(F)$ is a quadric cone in $\mathbb P^3.$
\end{description}
Moreover in all the cases above the pencil $|E|$ corresponds under the map
$\varphi_D$ to the system of generatrices of $\varphi_D(F).$
\end{theorem}

\section{Proof of Theorem \ref{propositionK3}}

We say that a curve $D$ is {\em nef} and {\em big} if $DC\geq 0$ for every curve $C$
and $D^2>0.$\newline
In order to prove Theorem \ref{propositionK3}, we show the following:

\begin{proposition}\label{K3}
Let $P$ be a smooth $K3$ surface with a reduced curve $B$ satisfying:
\begin{itemize}
    \item [{\rm (i)}] $B=B'+\sum_1^tA_i,$ $t\in \{ 9,\ldots,16\},$ where
    $B'$ is a nef and big curve with at most negligible
    singularities, the curves $A_i$ are disjoint $(-2)$-curves also
    disjoint from $B'$ and $B\equiv 2L,$ $L^2=-4,$
    for some $L\in {\rm Pic}(P).$
\end{itemize}
Then:
\begin{itemize}
    \item [{\rm a)}] Let $\pi:V\rightarrow P$ be a
    double cover with branch locus $B$ and $S$ be the smooth minimal
    model of $V$. Then $q(S)=0,$ $p_g(S)=1,$ $K_S^2=t-8$ and 
    the bicanonical map of $S$ is
    composed with the involution $i$ of $S$ induced by $\pi;$
    \item [{\rm b)}] If $t\geq 10,$ then $P$ contains a smooth curve
    $B_0'$ and $(-2)$-curves $A'_1,\ldots,A'_{t-1}$ such that
    $B_0'^2=B'^2-2$ and
    $B_0:=B_0'+\sum_1^{t-1}A'_i$ also satisfies condition {\rm (i)}.
\end{itemize}
\end{proposition}
{\bf Proof:}\\\\
a) Let $L\equiv\frac{1}{2}B$ be the line bundle which determines $\pi.$
From the double cover formulas (see {\it e.g.} \cite{BPV}) and the Riemann Roch theorem,
$$q(S)=h^1(P,\mathcal{O}_P(L)),$$
$$p_g(S)=1+h^0(P,\mathcal{O}_P(L)),$$
$$h^0(P,\mathcal{O}_P(L))+h^0(P,\mathcal{O}_P(-L))=h^1(P,\mathcal{O}_P(L)).$$
Since $2L-\sum A_i$ is nef and big, the Kawamata-Viehweg's vanishing Theorem
(see {\it e.g.} \cite[Corollary 5.12, c)]{EV}) implies $h^1(P,\mathcal O_P(-L))=0.$
Hence $$h^1(P,\mathcal{O}_P(L))=h^1(P,\mathcal{O}_P(K_P-L)=h^1(P,\mathcal O_P(-L)))=0$$
and then $q(S)=0$ and $p_g(S)=1.$
As $$h^0(P,\mathcal{O}_P(2K_P+L))=h^0(P,\mathcal{O}_P(L))=0,$$
the bicanonical map of $S$ is composed with $i$ (see \cite[Proposition 6.1]{CM}).

The $(-2)$-curves $A_1,\ldots,A_t$ give rise to $(-1)$-curves in $V,$ therefore
$$K_S^2=K_V^2+t=2(K_P+L)^2+t=2L^2+t=t-8.$$
\\
\\
b) Denote by $\xi\subset P$ the set of irreducible curves which do not
    intersect $B'$ and denote by $\xi_i,$ $i\geq 1,$ the connected components of $\xi.$
    Since $B'^2\geq 2,$ the Hodge-index Theorem implies
    that the intersection matrix of the components of $\xi_i$ is negative
    definite. Therefore, following \cite[Lemma I.2.12]{BPV}, the $\xi_i$'s have
    one of the five configurations: the support of $\textsf{A}_n,$ $\textsf{D}_n,$
    $\textsf{E}_6,$ $\textsf{E}_7$ or $\textsf{E}_8$ (see {\it e.g.} \cite[III.3]{BPV}
    for the description of these graphs).\\\\
    {\bf Claim 1:} {\em Each nodal curve $A_i$ can only be contained
    in a graph of type $\textsf{A}_{2n+1}$ or $\textsf{D}_n.$}\\\\
    {\em Proof}\ : Suppose that there exists an $A_i$ which is contained
    in a graph of type $\textsf{E}_6.$
    Denote the components of $\textsf{E}_6$ as in Figure \ref{E6}.

 \begin{figure}[htp]
\begin{center}
\begin{picture}(160,60)(0,-10)\label{picture}
\put(0,30){\line(1,0){160}}
\put(80,30){\line(0,-1){40}}
\put(0,30){\circle*{3.7}}
\put(40,30){\circle*{3.7}}
\put(80,30){\circle*{3.7}}
\put(120,30){\circle*{3.7}}
\put(160,30){\circle*{3.7}}
\put(80,-10){\circle*{3.7}}
\put(0,40){\makebox(0,0){$a_1$}}
\put(40,40){\makebox(0,0){$a_2$}}
\put(80,40){\makebox(0,0){$a_3$}}
\put(120,40){\makebox(0,0){$a_4$}}
\put(160,40){\makebox(0,0){$a_5$}}
\put(80,-20){\makebox(0,0){$a_6$}}
\end{picture}
\end{center}
\caption{$\textsf E_6$}
\label{E6}
\end{figure}

\noindent    If $A_i=a_3$ or $A_i=a_6,$ then $a_6B=a_6a_3=1$ or $a_3B=1,$
    contradicting $B\equiv 2L.$
    If $A_i=a_1$ or $A_i=a_2,$ then $a_2B=1$ or $a_1B=1,$ the same
    contradiction.
    By the same reason, $A_i\not=a_4$ and $A_i\not=a_5.$

    Analogously one can verify that each $A_i$ can not be in a
    graph of type $A_{2n},$ $E_7$ or $E_8.$ $\diamondsuit$
    \\

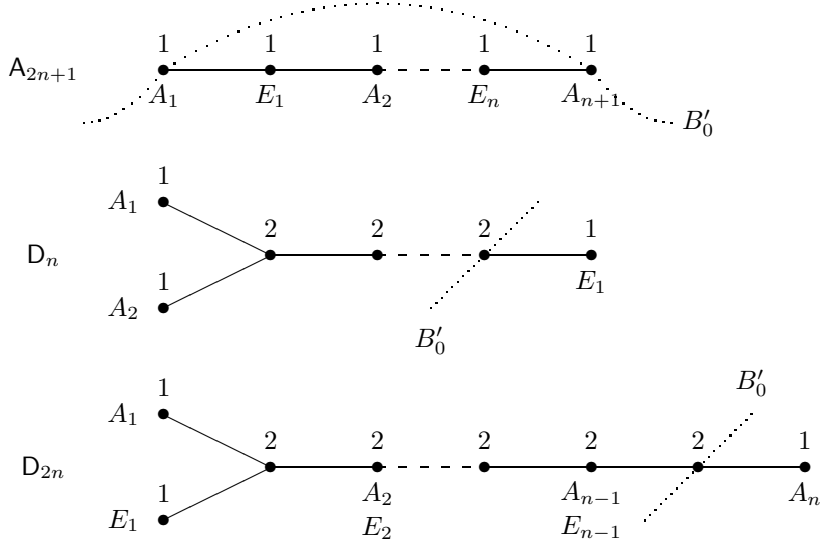
\begin{figure}[htp]
\begin{center}
\begin{picture}(250,200)(0,0)
\put(0,0){\line(2,1){40}}
\put(0,40){\line(2,-1){40}}
\put(40,20){\line(1,0){40}}
\multiput(80,20)(8,0){5}{\line(1,0){3}}
\put(120,20){\line(1,0){120}}
\put(0,0){\circle*{3.7}}
\put(0,40){\circle*{3.7}}
\put(40,20){\circle*{3.7}}
\put(80,20){\circle*{3.7}}
\put(120,20){\circle*{3.7}}
\put(160,20){\circle*{3.7}}
\put(200,20){\circle*{3.7}}
\put(240,20){\circle*{3.7}}
\put(0,10){\makebox(0,0){1}}
\put(0,50){\makebox(0,0){1}}
\put(40,30){\makebox(0,0){2}}
\put(80,30){\makebox(0,0){2}}
\put(120,30){\makebox(0,0){2}}
\put(160,30){\makebox(0,0){2}}
\put(200,30){\makebox(0,0){2}}
\put(240,30){\makebox(0,0){1}}
\put(80,10){\makebox(0,0){$A_2$}}
\put(80,-3){\makebox(0,0){$E_2$}}
\put(160,10){\makebox(0,0){$A_{n-1}$}}
\put(160,-3){\makebox(0,0){$E_{n-1}$}}
\put(240,10){\makebox(0,0){$A_{n}$}}
\put(-15,0){\makebox(0,0){$E_1$}}
\put(-15,40){\makebox(0,0){$A_1$}}
\put(-45,20){\makebox(0,0){$\textsf D_{2n}$}}
\put(-45,100){\makebox(0,0){$\textsf D_{n}$}}
\put(0,80){\line(2,1){40}}
\put(0,120){\line(2,-1){40}}
\put(40,100){\line(1,0){40}}
\multiput(80,100)(8,0){5}{\line(1,0){3}}
\put(120,100){\line(1,0){40}}
\put(0,80){\circle*{3.7}}
\put(0,120){\circle*{3.7}}
\put(40,100){\circle*{3.7}}
\put(80,100){\circle*{3.7}}
\put(120,100){\circle*{3.7}}
\put(160,100){\circle*{3.7}}
\put(0,90){\makebox(0,0){1}}
\put(0,130){\makebox(0,0){1}}
\put(40,110){\makebox(0,0){2}}
\put(80,110){\makebox(0,0){2}}
\put(120,110){\makebox(0,0){2}}
\put(160,110){\makebox(0,0){1}}
\put(160,90){\makebox(0,0){$E_1$}}
\put(-15,80){\makebox(0,0){$A_2$}}
\put(-15,120){\makebox(0,0){$A_1$}}
\qbezier[20](180,0)(200,20)(220,40)
\qbezier[20](100,80)(120,100)(140,120)
\put(220,52){\makebox(0,0){$B_0'$}}
\put(100,68){\makebox(0,0){$B_0'$}}
\put(0,170){\line(1,0){80}}
\put(120,170){\line(1,0){40}}
\multiput(80,170)(8,0){5}{\line(1,0){3}}
\put(0,170){\circle*{3.7}}
\put(40,170){\circle*{3.7}}
\put(80,170){\circle*{3.7}}
\put(120,170){\circle*{3.7}}
\put(160,170){\circle*{3.7}}
\put(0,180){\makebox(0,0){1}}
\put(40,180){\makebox(0,0){1}}
\put(80,180){\makebox(0,0){1}}
\put(120,180){\makebox(0,0){1}}
\put(160,180){\makebox(0,0){1}}
\put(0,160){\makebox(0,0){$A_1$}}
\put(40,160){\makebox(0,0){$E_1$}}
\put(80,160){\makebox(0,0){$A_2$}}
\put(120,160){\makebox(0,0){$E_n$}}
\put(160,160){\makebox(0,0){$A_{n+1}$}}
\put(-45,170){\makebox(0,0){$\textsf A_{2n+1}$}}
\qbezier[40](0,170)(80,220)(160,170)
\qbezier[10](190,150)(175,150)(160,170)
\qbezier[10](-30,150)(-15,150)(0,170)
\put(200,150){\makebox(0,0){$B_0'$}}
\end{picture}
\end{center}
\caption{\em The numbers represent the multiplicity and the doted curve
represent a general element $B_0'$ in $|B'-\textsf{G}|.$}
\label{ADD}
\end{figure}

    The possible configurations for the curves $A_i$ in the graphs are
    shown in Figure \ref{ADD}.
    Fix one of the curves $A_i$ and denote by $\textsf{G}$ the graph containing it.\\

\noindent
    {\bf Claim 2:} {\em We can choose $A_i$ such that the linear
    system $|B'-\textsf{G}|$ has no fixed components (and thus no base points,
    from {\rm \cite[Theorem 3.1]{St}}).}\\\\
    {\em Proof}\ : Denote by $\varphi_{|B'|}$ the map given by the linear system
    $|B'|.$ We know that $\varphi_{|B'|}$ is birational or it is of
    degree 2 (see \cite[Section 4]{St}). If $\varphi_{|B'|}$ is birational
    or the point $\varphi_{|B'|}(\textsf{G})$ is a smooth point of
    $\varphi_{|B'|}(P),$ the result is clear, since $|B'-\textsf{G}|$ is
    the pullback of the linear system of the hyperplanes containing
    $\varphi_{|B'|}(\textsf{G})$ and
    $\varphi_{|B'|}^*(\varphi_{|B'|}(\textsf{G}))=\textsf{G}$
    (see \cite[Theorems III 7.1 and 7.3]{BPV}).

    Suppose now that $\varphi_{|B'|}$ is non-birational and that
    $\varphi_{|B'|}(\textsf{G})$ is a singular point of
    $\varphi_{|B'|}(P).$ Then $B'$ is linearly equivalent to a curve
    with one of the configurations described in Theorem \ref{SD}.
    Except for the last configuration, $\textsf{G}$ contains at
    most two $(-2)$-curves. But $t\ge 9,$ thus in these cases there exists other
    graph $\textsf{G}'$ containing a curve $A_j$ such that
    $\varphi_{|B'|}(\textsf{G}')$ is a non-singular point of $\varphi_{|B'|}(P)$
    (notice that Theorem \ref{SD} implies that $\varphi_{|B'|}(P)$ contains
    only one singular point).

    So we can suppose that $B'$ is equivalent to a curve with
    a configuration as in Theorem \ref{SD}, (iii), b). None of the curves
    $\Gamma_0,\ldots,\Gamma_N$ can be
    one of the curves $A_j.$ For this note that: if $\Gamma_0=A_j,$ then
    $EB=E\left(B'+\sum A_i\right)=2+E\Gamma_0=3\not\equiv 0\ ({\rm mod\ 2});$
    if $\Gamma_1=A_j,$ then $\Gamma_0 B=\Gamma_0\Gamma_1=1\not\equiv 0\
    ({\rm mod\ 2});$ etc. Again this configuration can contain at most two
    curves $A_j,$ the components $\Gamma_{N+1},$ $\Gamma_{N+2}.$ $\diamondsuit$
    \\

    Let $B_0'$ be a smooth curve in $|B'-\textsf{G}|.$
    If $\textsf{G}$ is an $\textsf{A}_{2n+1}$ graph, then,
    using the notation of Figure \ref{ADD},
    $$\left( B_0'+\sum_1^nE_i\right) +\sum_{n+2}^tA_i\equiv \left( B'-\sum_1^{n+1}A_i\right) +\sum_{n+2}^tA_i\equiv$$
    $$\equiv B'+\sum_1^tA_i-2\sum_1^{n+1}A_i\equiv 0\ {\rm (mod\ 2)}.$$
    Therefore the curve $$B_0:=B_0'+\sum_1^nE_i+\sum_{n+2}^tA_i$$ satisfies condition (i).

    The case where $\textsf{G}$ is a $\textsf{D}_{m}$ graph is analogous.
\Qed
{\bf Proof of Theorem \ref{propositionK3}\ :}
Let $V\rightarrow S$ be the blow-up at the isolated fixed points of the involution $i$ and $W$ be the minimal resolution
of $S/i.$ We have a commutative diagram
$$
\begin{CD}\ V@> >>S\\ @V\pi VV  @VV  V\\ W@> >> S/i \ .
\end{CD}
$$
The branch locus of $\pi$ is a smooth curve $B=B'+\sum_1^t A_i,$ where the curves
$A_i$ are $(-2)$-curves which contract to the nodes of $S/i.$
Let $P$ be the minimal model of $W$ and $\overline B\subset P$ be the projection of $B.$
Let $L\equiv\frac{1}{2}B$ be the line bundle which determines $\pi.$

First we verify that $\overline B$ satisfies condition (i) of Proposition \ref{K3}:
from \cite[Proposition 6.1]{CM}, $\chi(\mathcal O_W)-\chi(\mathcal O_S)=K_W(K_W+L),$
hence $K_W(K_W+L)=0,$ which implies that $\overline B$ has at most negligible singularities; 
now from \cite[Theorem 5.2]{Mo} we get $K_S^2=\frac{1}{2}\overline{B'}^2$ and 
$1=p_g(S)=\frac{1}{4}(K_S^2-t)+3,$ thus $t=K_S^2+8$ and 
$\overline B^2=\overline{B'}^2-2t=2K_S^2-2t=-16,$ which gives $(\overline B/2)^2=-4$ and $\overline{B'}^2\geq 2;$
finally $\overline{B'}$ is nef because, on a $K3$ surface, an irreducible curve with negative self intersection
must be a $(-2)$-curve.

Now using Proposition \ref{K3}, b) and a) we obtain statement b).
In particular we get also that $P$
contains a curve $B_0'$ and $(-2)$-curves $A'_i,$ $i=1,\ldots,9,$ such that
$B_0:=B_0'+\sum_1^9A'_i$ is smooth and divisible by 2 in the Picard group.
Moreover, the complete linear system $|B_0'|$
has no fixed component nor base points and $B_0'^2=2.$
Therefore, from \cite{St}, $|B_0'|$ defines a
generically finite degree $2$ morphism $$\varphi:=\varphi_{|B_0'|}:P\rightarrow\mathbb P^2.$$
Since $g(B'_0)=2,$ this map is ramified over a sextic curve $\beta$.
The singularities of $\beta$ are negligible because $P$ is a $K3$ surface.

We claim that $\beta$ is the union of two cubics.
Let $p_i\in\beta$ be the singular point corresponding to $A'_i,$ $i=1,\ldots,9.$
Notice that the $p_i$'s are possibly infinitely near.
Let $C\subset\mathbb P^2$ be a cubic curve passing through $p_i,$ $i=1,\ldots,9.$
As $C+\varphi_*(B_0')$ is a plane quartic, we have
$$\left(\varphi^*(C)-\sum_1^9A'_i\right)+B_0'+\sum_1^9A'_i\equiv \varphi^*
(C+\varphi_*(B_0'))\equiv 0\ {\rm (mod\ 2)},$$
hence also $\varphi^*(C)-\sum_1^9A'_i\equiv 0\ {\rm (mod\ 2)},$ {\it i.e.}
there exists a divisor $J$ such that $$2J\equiv\varphi^*(C)-\sum_1^9A'_i.$$
Since $P$ is a $K_3$ surface,
the Riemann Roch theorem implies that $J$ is effective. From $JA'_i=1,$ $i=1,\ldots,9,$
we obtain that the plane curve $\varphi_*(J)$ passes with multiplicity 1 through the
nine singular points $p_i$ of $\beta.$
This immediately implies that $\varphi_*(J)$ is not a line nor a conic,
because $\beta$ is a reduced sextic.
Therefore $\varphi_*(J)$ is a reduced cubic. So $\varphi_*(J)\equiv C$ and then
$$\varphi^*\left(\varphi_*(J)\right)\equiv 2J+\sum_1^9A_i'.$$
This implies that $\varphi_*(J)$ is contained in the branch locus $\beta,$
which finishes the proof of a).\Qed

\section{Examples}

Todorov gave examples of surfaces $S$ with bicanonical image $\phi_2(S)$
birational to a Kummer surface having only ordinary double points as
singularities. The next sections contain an example with $\phi_2(S)$
non-birational to a Kummer surface and an example with $\phi_2(S)$ having
an $\textsf A_{17}$ double point.

\subsection{$S/i$ non-birational to a Kummer surface}\label{explononkummer}

Here we construct smooth surfaces $S$ of general type with $K^2=2,3,$
$p_g=1$ and $q=0$ having bicanonical map of degree 2 onto a $K3$ surface
which is not birational to a Kummer surface.\\

It is known since \cite{Hu} that there exist special sets of 6 nodes, called
Weber hexads, in the Kummer surface $Q\in\mathbb P^3$ such that the surface
which is the blow-up of $Q$ at these nodes can be embedded in $\mathbb P^3$
as a quartic with 10 nodes. This quartic is the Hessian of a smooth cubic
surface. 

The space of all smooth cubic surfaces has dimension 4 while the
space of Kummer surfaces has dimension 3. Thus it is natural to ask
if there exist Hessian "non-Kummer" surfaces, {\it i.e.} which are not the
embedding of a Kummer surface blown-up at 6 points. This is studied
in \cite{Ro}, where the existence of "non-Kummer" quartic Hessians
$H$ in $\mathbb P^3$ is shown. These are surfaces with 10 nodes $a_i$ such that
the projection from one node $a_1$ to $\mathbb P^2$ is a generically $2:1$
cover of $\mathbb P^2$ with branch locus $\alpha_1+\alpha_2$ satisfying:
$\alpha_1,$ $\alpha_2$ are smooth cubics tangent to a nondegenerate
conic $C$ at 3 distinct points. We use this in the following construction.\\

Let $\alpha_1,$ $\alpha_2$ and $C$ be as above.
Take the morphism $\pi:W\rightarrow\mathbb P^2$
given by the canonical resolution of the double cover of $\mathbb P^2$
with branch locus $\alpha_1+\alpha_2.$
The strict transform of $C$ gives
rise to the union of two disjoint $(-2)$-curves $A_1,A_2\subset W$
(one of these correspond to the node $a_1$ from which we have projected).

Let $T\in\mathbb P^2$ be a general line. Let $A_3,\ldots,A_{11}\subset W$ be
the disjoint $(-2)$-curves contained in $\pi^*(\alpha_1+\alpha_2).$
We have $\pi^*(T+\alpha_1)\equiv 0\ {\rm (mod\ 2)},$
hence, since $\alpha_1$ is in the branch locus, also
$$\pi^*(T)+\sum_3^{11} A_i\equiv 0\ {\rm (mod\ 2)}.$$
The linear systems $|\pi^*(T)+A_2|$ and $|\pi^*(T)+A_1+A_2|$ have no fixed components
nor base points (see \cite[(2.7.3) and Corollary 3.2]{St}).
The surface $S$ is the minimal model of the double cover of $W$ ramified over a general element in
$$|\pi^*(T)+A_2|+\sum_2^{11}A_i \ \ \ \ {\rm or}\ \ \ \ |\pi^*(T)+A_1+A_2|+\sum_1^{11}A_i.$$

\subsection{$\phi_2(S)$ with $A_{17}$ and $A_1$ singularities}

This section contains a brief description of a construction of a surface $S$
of general type having bicanonical image $\phi_2(S)\subset\mathbb P^3$
a quartic $K3$ surface with $A_{17}$ and $A_1$ singularities. 
I ommit the details, which were verified using the 
{\em Computational Algebra System Magma}.\\

Let $C_1$ be a nodal cubic, $p$ an inflection point of $C_1$ and $T$
the tangent line to $C_1$ at $p.$ The pencil generated by $C_1$ and
$3T$ contains another nodal cubic $C_2,$ smooth at $p.$ The curves
$C_1$ and $C_2$ intersect at $p$ with multiplicity 9.

Let $\rho:X\rightarrow\mathbb P^2$ be the resolution of $C_1+C_2$ and $\pi:W\rightarrow X$
be the double cover with branch locus the strict transform of $C_1+C_2.$
Denote by $\overline l$ the line containing the nodes of
$C_1$ and $C_2$ and by $l\subset W$ the pullback of the strict transform of $\overline l.$
The map given by $|(\rho\circ\pi)^*(\overline l)+l|$ is birational onto a
quartic $Q$ in $\mathbb P^3$ with an $A_1$ and $A_{17}$ singularities
(notice that $l$ is a $(-2)$-curve and $((\rho\circ\pi)^*(\overline l)+l)l=0$).

Let $B'\in |(\rho\circ\pi)^*(\overline l)+l|$ be a smooth element and $A_1,\ldots,A_9$ be
the disjoint $(-2)$-curves contained in $(\rho\circ\pi)^*(p).$
Let $S$ be the minimal model of the double cover of
$W$ with branch locus $B'+\sum_1^9A_i+l.$
The surface $Q$ is the image of the bicanonical map of
$S$ and $p_g(S)=1,$ $q(S)=0,$ $K_S^2=2.$
\bibliography{ReferencesRito2}

\bigskip
\bigskip

\noindent Carlos Rito
\\ Departamento de Matem\' atica
\\ Universidade de Tr\' as-os-Montes e Alto Douro
\\ 5000-911 Vila Real
\\ Portugal
\\\\
\noindent {\it e-mail:} crito@utad.pt

\end{document}